# The density of ramified primes in semisimple $p$-adic Galois representations

Chandrashekhar Khare and C. S. Rajan

## 1 Introduction

Let $L$ be a number field. Consider a continuous, semisimple $p$-adic Galois Galois representation

$$\rho : G_L \to GL_m(K)$$

of the absolute Galois group $G_L$ of $L$, and with $K$ a finite extension of $\mathbf{Q}_p$. In [R] in the case when $n = 2$ and $L = \mathbf{Q}$ examples of such representations were constructed that were ramified at infinitely many primes (see also the last section of [KR]: we call such representations *infinitely ramified*), and which had open image and determinant $\varepsilon$ the $p$-adic cyclotomic character. One may ask if in these examples of [R] the set of ramified primes is of small density.

**Theorem 1** *Let $\rho : G_L \to GL_m(K)$ be a continuous, semisimple representation. Then the set of primes $S_\rho$ that ramify in $\rho$ is of density zero.*

The semisimplicity assumption is crucial as a construction using Kummer theory (see exercise on III-12 of [S]) gives examples of continuous, reducible but indecomposable representations $\rho : G_L \to GL_2(\mathbf{Q}_p)$ that are ramified at *all primes.* Note that as in [R], infinitely ramified representations, though not being *motivic* themselves, do arise as limits of motivic $p$-adic representations.

The results of this paper are used in [K] to study the fields of rationality in a converging sequence of *algebraic*, semisimple $p$-adic Galois representations. After Theorem 1, we know that the set of primes which are unramified in a continuous, semisimple representation $\rho : G_L \to GL_m(K)$ is of density one. Hence many of the results (e.g., the strong multiplicity one results of [Ra]), available in the *classical* case when $\rho$ is assumed to be finitely ramified, extend to this more general situation. After Theorem 1 it also makes good sense to talk of compatible systems of continuous, semisimple Galois representations in the sense of [S], without imposing the condition that these be finitely ramified. We raise the following question:

**Question 1** *Given two compatible continuous, semisimple representations $\rho : G_L \to GL_m(\mathbf{Q}_\ell)$ and $\rho' : G_L \to GL_m(\mathbf{Q}_{\ell'})$ with $\ell \neq \ell'$, then is the set of primes at which either of $\rho$ or $\rho'$ ramifies finite?*

## 2 Proof of theorem

### 2.1

Let $\rho$ be as in the theorem. As $\rho$ is continuous, we can regard $\rho$ as taking values in $GL_m(\mathcal{O})$ where $\mathcal{O}$ is the ring of integers of $K$. We denote the maximal ideal of $\mathcal{O}_K$ by $\mathfrak{m}$, and by $\rho_n$ the reduction mod $\mathfrak{m}^n$ of $\rho$.

We define $c_{\rho,n}$ to be the upper density of the set $S_{\rho,n}$ of primes $q$ of $L$ that

1. lie above primes which split in $L/\mathbf{Q}$ (this assumption is merely for notational convenience: we denote abusively the prime of $\mathbf{Q}$ lying below it by $q$),

2. are unramified in $\rho_1$ and $\neq p$,

3. $\rho_n|_{D_q}$ is unramified, but there exists a "lift" of $\rho_n|_{D_q}$, with $D_q$ the decomposition group at $q$ corresponding to a place above $q$ in $\overline{\mathbf{Q}}$, to a representation $\tilde{\rho}_q$ of $D_q$ to $GL_m(K)$ that is ramified at $q$: by a lift we mean some conjugate of $\tilde{\rho}_q$ which reduces mod $\mathfrak{m}^n$ to $\rho_n|_{D_q}$. Note that by 2, any such lift $\tilde{\rho}_q$ factors through $G_q$, the quotient of $D_q$ that is the Galois group of the maximal tamely ramified extension of $\mathbf{Q}_q$.

**Proposition 1** *Given any $\varepsilon > 0$, there is an integer $N_\varepsilon$ such that $c_{\rho,n} < \varepsilon$ for $n > N_\varepsilon$.*

We claim that the proposition implies Theorem 1. To see this first observe that the primes of $L$ that do not lie above primes of $\mathbf{Q}$ which split in the extension $L/\mathbf{Q}$ are of density 0. To prove the theorem it is enough to show that given any $\varepsilon > 0$, the upper density of the set $S_\rho$ of ramified primes for $\rho$ is $< \varepsilon$. Consider the $N_\varepsilon$ that the proposition provides. Further note that in $\rho_{N_\varepsilon}$ only finitely many primes ramify. From this it readily follows that the upper density of $S_\rho$ is $< \varepsilon$. Hence the theorem. Thus it only remains to prove the proposition.

## 2.2 Tame inertia

The proposition relies on the structure of the Galois group $G_q$ of the maximal tamely ramified extension of $\mathbf{Q}_q$. This is used to calculate the densities $c_{\rho,n}$ for large enough $n$. The concept of largeness of $n$ for our purposes is independent of the representation $\rho$ and the prime $q$, and depends only on $K$ and the dimension of the representation. Roughly the idea of the proof of the proposition is that for these large $n$, only *semistable* (i.e., image of inertia is unipotent) lifts intervene in the calculation of the $c_{\rho,n}$, and the conjugacy classes in the image of $\rho_n$ of the Frobenii of the primes in $S_{\rho,n}$ can be seen to lie in the $\mathcal{O}/\mathsf{m}^n$ valued points of a analytically defined subset of $\text{im}(\rho)$ of smaller dimension. We flesh out this idea below. We will implicitly use the fact that though one cannot speak of eigenvalues of an element of $GL_m(A)$, for a general ring $A$, its characteristic polynomial makes good sense.

The group $G_q$ is topologically generated by two elements $\sigma_q$ and $\tau_q$ that satisfy the relation

$$\sigma_q \tau_q \sigma_q^{-1} = \tau_q^q, \tag{1}$$

and such that $\sigma_q$ induces the Frobenius on residue fields and $\tau_q$ (topologically) generates the tame inertia subgroup.

## 2.3 Reduction to the semistable case

**Lemma 1** *Let $\theta : G_q \to GL_m(K)$ be any continuous representation. Then the roots of the characteristic polynomial of $\theta(\tau_q)$ are roots of unity. Further the order of these roots is bounded by a constant depending only on $K$.*

**Proof:** Using Krasner's lemma, we know that there are only finitely many degree $m$ extensions of $K$. Let $K'$ be the finite extension of $K$ that is the compositum of all the degree $m$ extensions of $K$. By extending scalars to $K'$, we can assume that $\theta(\tau_q)$ is upper triangular. Let $\theta_1, \cdots, \theta_m$ be the diagonal entries. Using equation (1), we deduce that

$$\{\theta_1, \cdots, \theta_m\} = \{\theta_1^q, \cdots, \theta_m^q\}.$$

From this it follows that the $\theta_i$ are roots of unity (of order dividing $q^m - 1$). Hence the last statement of the lemma follows from the fact that there are only finitely many roots of unity in $K'$.

**Corollary 1** *Let $\theta : G_q \to GL_m(K)$ be any continuous representation. Assume that the characteristic polynomial of $\theta(\tau_q)$ is not equal to $(x-1)^m$. Then there exists an integer $N(m, K)$ depending only on $m$ and $K$, such that the reduction modulo $\mathsf{m}^{N(m,K)}$ of any conjugate of $\theta$ into $GL_m(\mathcal{O})$ is ramified.*

**Proof:** Choose $N(m, K)$ such that if $\zeta \in K'^*$ is a root of unity satisfying $(\zeta - 1)^m \equiv 0 \pmod{\mathsf{m}^{N(m,K)}}$, then $\zeta = 1$, for $K'$ as in the proof of lemma above. Then the corollary follows by considering reductions of characteristic polynomials.

**Corollary 2** *In a continuous, semisimple representation of $\rho : G_L \to GL_m(K)$, the set of primes $q$ as above for which $\rho(\tau_q)$ is not unipotent, is finite.*

**Corollary 3** *Any continuous, semisimple abelian representation of $G_L \to GL_m(K)$ is finitely ramified.*

## 2.4 The $GL_2$ case

At this point for the sake of exposition, we briefly indicate the proof of the theorem when $m = 2$, and $\rho(G_L)$ is open in $GL_2(K)$ with determinant $\varepsilon$ the $p$-adic cyclotomic character. (Note that in the case when the Lie algebra of $\rho(G_L)$ is abelian, the ramification set is finite by Corollary 3.)

Consider $S_{\rho,n}$ for $n > N(2, K)$ and let $q \in S_{\rho,n}$. Let $\tilde{\rho}_q$ be any lift of $\rho_n|_{D_q}$ to $GL_2(K)$ that is ramified at $q$. By the above considerations, it follows that $\tilde{\rho}_q(\tau_q)$ is unipotent, which we can assume to be upper triangular. Since $\tilde{\rho}_q(\sigma_q)$ normalises $\tilde{\rho}_q(\tau_q)$, we can assume that

- $\tilde{\rho}_q(\tau_q)$ of the form $\begin{pmatrix} 1 & * \\ 0 & 1 \end{pmatrix}$ and $\tilde{\rho}_q(\tau_q)$ non-trivial

- $\tilde{\rho}_q(\sigma_q)$ is of the form $\begin{pmatrix} \alpha & * \\ 0 & \beta \end{pmatrix}$.

Observe that $\alpha \neq \beta$ because of the relation (1). Thus we can further assume by conjugating by an element of the form $\begin{pmatrix} 1 & 1 \\ 0 & y \end{pmatrix} \in GL_2(K)$ that $\tilde{\rho}_q(\sigma_q)$ is of the form $\begin{pmatrix} \alpha & 0 \\ 0 & \beta \end{pmatrix}$. Then we see from the equation (1), that $\alpha \beta^{-1} = q$.

Consider the invariant functions tr and det defined on the space of conjugacy classes of $GL_2(\mathcal{O})$ or $GL_2(\mathcal{O}/\mathsf{m}^n)$ given by the trace and determinant functions. We see from our work above that primes $q \in S_{\rho,n}$ for $n > N(2, K)$ are such that the conjugacy classes $\rho_n(\mathrm{Frob}_q)$ satisfy the relation

$$\mathrm{tr}^2 = (1 + \det)^2.$$

From this we conclude using the fact that the image of $\rho$ is open in $GL_2(K)$, Cebotarev density theorem and the second paragraph on page 586 of [S1] that $c_{n,\rho} \to 0$ as $n \to \infty$. The proposition follows in this case and the proof of the theorem is complete in the special case of open image in $GL_2(K)$ with determinant $\varepsilon$.

## 2.5 The general case

We reduce the general situation to the case when $\mathrm{im}(\rho)$ is a semisimple $p$-adic Lie group contained in $GL_M(\mathbf{Q}_p)$ for some $M$. Firstly by Weil restriction of scalars we may assume $K = \mathbf{Q}_p$ (with possibly a different $m$). Let $G$ be the Zariski closure of the image of $\rho$. Since $\rho$ is semisimple, $G$ is reductive and let $Z$ be the centre of the connected component of $G$. Let $\rho_s : G_L \to (G/Z)(\mathbf{Q}_p)$ be the corresponding representation. Because of Corollary 2 we see that the ramification set of $\rho$ and $\rho_s$ differ by a finite set, and thus we can work with $\rho_s$. Now embed $G/Z$ into $GL_M/\mathbf{Q}_p$ for some $M$. Thus we have reduced to the case when $\mathrm{im}(\rho)$ is a semisimple $p$-adic Lie group contained in $GL_M(\mathbf{Q}_p)$ for some $M$.

We look at $S_{\rho,n}$ for $n > N(M, \mathbf{Q}_p)$ and let $q \in S_{\rho,n}$. Let $\tilde{\rho}_q$ be any lift of $\rho_n|_{D_q}$ to $GL_M(K)$ that is ramified at $q$. By Corollary 1, we can assume that $\tilde{\rho}_q(\tau_q)$ is unipotent (*and non-trivial*), which we can further assume to be upper triangular.

Consider the canonical filtration of $\tilde{\rho}_q(\tau_q)$ acting on the vector space $\mathbf{Q}_p^M$, with the dimension of the corresponding graded components $m_1, \cdots, m_i$. By conjugating by an element in the Levi, over a finite extension of $\mathbf{Q}_p$, of the corresponding parabolic subgroup defined by $\tilde{\rho}_q(\tau_q)$ (of the form $GL_{m_1} \times \cdots \times GL_{m_i}$), we can assume that $\tilde{\rho}_q(G_q)$ is upper triangular.

**Lemma 2** *If $f_q(x)$ is the characteristic polynomial of $\tilde{\rho}_q(\sigma_q)$, then $f_q(x)$ and $f_q(qx)$ have a common root.*

**Proof:** Let $U$ be the subgroup of unipotent upper triangular matrices of $SL_M(\mathbf{Q}_p)$, and let
$$U = U^0 \supset U^1 \supset \cdots \supset 1$$
be the descending central filtration. Let $i$ be the smallest integer such that $\tilde{\rho}_q(\tau_q) \notin U^{i+1}$. By looking at the conjugation action of $\tilde{\rho}_q(\sigma_q)$ on $U^i/U^{i+1}$, and using the relation (1), it follows that there are two eigenvalues $\alpha_q$, $\beta_q$ of $\tilde{\rho}_q(\sigma_q)$ such that $\alpha_q \beta_q^{-1} = q$. Hence the lemma.

Consider $\rho' = \rho \oplus \varepsilon : G_L \to GL_M(\mathbf{Q}_p) \times GL_1(\mathbf{Q}_p)$. Let $G' = G \times GL_1$. Choose an integral model for $\rho'$, i.e., $\rho'(G_L) \subset GL_M(\mathbf{Z}_p) \times GL_1(\mathbf{Z}_p)$, induced by the chosen integral model of $\rho$, and denote by $\rho'_n$ its reduction mod $\mathsf{m}^n$. We normalise the isomorphism of class field theory so that a uniformiser is sent to the arithmetic Frobenius (so $\varepsilon(\mathrm{Frob}_q) = q$).

Let
$$(A, b) \in GL_M(\mathbf{Q}_p) \times GL_1(\mathbf{Q}_p)$$
and let $f(x)$ be the characteristic polynomial of $A$. Let $F$ be the invariant polynomial function on $GL_M \times GL_1$ with $\mathbf{Z}_p$-coefficients defined by the resultant of the two polynomials $f(x)$ and $f(bx)$.

By choosing $b$ different from the ratios of eigenvalues of an element of $G$, we deduce that no connected component of $G'$ is contained inside the variety $F = 0$. Thus we see that $\{F = 0\} \cap G'$ is a subvariety of smaller dimension than the dimension of $G'$.

**Lemma 3** *$\rho'(G_L)$ is an open subgroup of $G'(\mathbf{Q}_p) = G(\mathbf{Q}_p) \times GL_1(\mathbf{Q}_p)$.*

**Proof:** Since $\mathrm{im}(\rho)$ is a semisimple $p$-adic group, we deduce from Chevalley's theorem (Corollary 7.9 of [Bo]) that $\mathrm{im}(\rho)$ is open in $G(\mathbf{Q}_p)$. From this we further deduce that the commutator subgroup of $\mathrm{im}(\rho)$ is of finite index in $\mathrm{im}(\rho)$. Thus the intersection of the fixed fields of the kernel of $\rho$ and $\varepsilon$ is a finite extension of $\mathbf{Q}$. Certainly $\text{im}(\varepsilon)$ is open in $\mathbf{Q}_p^*$ and hence the lemma follows.

From the openness of $\text{im}(\rho')$, we see that

$$\lim_{n\to\infty} \frac{|\text{im}(\rho_n')|}{p^{nd}}$$

is a non-zero positive constant, where $d$ is the dimension of $G'$.

On the other hand using the notation and results of Section 3 of [S1], if we denote by $\tilde{Y}_n$ the elements $x \in \text{im}(\rho_n')$ that satisfy $F(x) \equiv 0 \pmod{p^n}$, then from the second paragraph on page 586 of [S1] it follows that $|\tilde{Y}_n| = O(p^{n(d-\delta)})$ where $\delta$ is a positive constant. By Lemma 2 we see that $\rho_n'(\text{Frob}_q) \in \tilde{Y}_n$ for $q \in S_{\rho,n}$. Then applying the Cebotarev density theorem we conclude that

$$c_{\rho,n} \leq \frac{|\tilde{Y}_n|}{|\text{im}(\rho_n')|},$$

and hence $c_{n,\rho} \to 0$ as $n \to \infty$. This finishes the proof of Proposition 1 and hence that of Theorem 1.

**Remarks:**

- To prove the theorem instead of defining $S_{\rho,n}$ the way we did, we could have have worked with the smaller subset consisting of primes that are unramified in $\rho_n$, but ramified in the $\rho$ of Theorem 1. In the notation of page 586 of [S1] we would then be working with $Y_n$ rather than $\tilde{Y}_n$ as above. By Theorem 8 of Section 3 of [S1] we will obtain a better estimate $c_{\rho,n} = O(p^{-n})$. This may be useful to get more precise quantitative versions of Theorem 1. We have defined $S_{\rho,n}$ the way we have for its use in [K].

- An analog of Theorem 1 is valid for function fields of characteristic $\ell \neq p$, and the same proof works. On the other hand for function fields of characteristic $p$, Theorem 1 is false, and in this case there are examples of semisimple $p$-adic Galois representations ramified at *all places*. It is easy to construct such examples using the fact that the Galois group in this case has $p$-cohomological dimension $\leq 1$.

*Address of the authors*: School of Mathematics, TIFR, Homi Bhabha Road, Mumbai 400 005, INDIA. e-mail addresses: shekhar@math.tifr.res.in, rajan@math.tifr.res.in